\documentclass[11pt, twoside]{amsart}
\usepackage{amsmath, amsthm, amscd, amsfonts, amssymb, graphicx, color}
\usepackage[bookmarksnumbered, plainpages, backref]{hyperref}

\begin{document}

\title{Hilbert Matrix and Difference Operator of Order $m$}

\author{Murat Kiri\c{s}ci*, Harun Polat}

\address{[Murat Kiri\c{s}ci] Department of Mathematical Education, Hasan Ali Y\"{u}cel Education Faculty,
Istanbul University, Vefa, 34470, Fatih, Istanbul, Turkey \vskip 0.1cm }
\email{mkirisci@hotmail.com, murat.kirisci@istanbul.edu.tr}

\address{[Harun Polat] Mu\c{s} Alparslan University Art and Science Faculty, Department of
Mathematics, 49100 Mu\c{s},Turkey \vskip 0.1cm }
\email{h.polat@alparslan.edu.tr}

\thanks{*Corresponding author.}

\begin{abstract}

In this paper, firstly, some applications of Hilbert matrix in image processing and cryptology are mentioned and an algorithm
related to the Hilbert view of a digital image is given. In section 2, the new matrix domains are constructed and some properties are investigated.
Furthermore, dual spaces of new matrix domains are computed and matrix transformations are characterized.
In last section, examples of transformations of new spaces are given.
\end{abstract}

\maketitle

\pagestyle{plain} \makeatletter
\theoremstyle{plain}
\newtheorem{thm}{Theorem}[section]
\numberwithin{equation}{section}
\numberwithin{figure}{section}  
\theoremstyle{plain}
\newtheorem{pr}[thm]{Proposition}
\theoremstyle{plain}
\newtheorem{exmp}[thm]{Example}
\theoremstyle{plain}
\newtheorem{cor}[thm]{Corollary} 
\theoremstyle{plain}
\newtheorem{defin}[thm]{Definition}
\theoremstyle{plain}
\newtheorem{lem}[thm]{Lemma} 
\theoremstyle{plain}
\newtheorem{rem}[thm]{Remark}
\numberwithin{equation}{section}

\section{Introduction}

\subsection{Hilbert Matrix and Applications}

We consider the matrices $H$  and $H_{n}$ as follows:

\begin{eqnarray*}
H=\begin{bmatrix}
1 & 1/2 & 1/3 & 1/4 & \cdots \\
1/2 & 1/3 & 1/4 & \cdots \\
1/3 & 1/4 & \cdots \\
1/4 & \cdots \\
\vdots & \vdots
\end{bmatrix} \quad \textrm{ and } \quad H_{n}=\begin{bmatrix}
1 & 1/2 & 1/3 & 1/4 & \cdots & 1/n\\
1/2 & 1/3 & 1/4 & \cdots \\
1/3 & 1/4 & \cdots \\
\vdots & \vdots \\
1/n & \cdots 1/(2n-1)
\end{bmatrix}
\end{eqnarray*}

It is well known that these matrices are called the infinite Hilbert matrix and the $n \times n$ Hilbert matrix, respectively.
A famous inequality of Hilbert (\cite{HLP}, Section 9) asserts that the matrix $H$ determines a bounded linear operator on the Hilbert space of square summable complex sequences. Also, $n \times n$ Hilbert matrices are well known examples of extremely of ill-conditioned matrices.\\

Frequently, Hilbert matrices are used both mathematics and computational sciences. For examples, in image processing, Hilbert matrices
are commonly used. Any $2D$ array of natural numbers in the range $[0,n]$ for all $n\in \mathbb{N}$ can be viewed as a greyscale digital image.\\

We take the Hilbert matrix $H_{n}$($n \times n$ matrix). If we use the Mathematica, then we can write
\begin{eqnarray*}
hilbert=HilbertMatrix[5]//MatrixForm
\end{eqnarray*}
and we can obtain
\begin{eqnarray*}
H=\begin{bmatrix}
1 & 1/2 & 1/3 & 1/4 & 1/5 \\
1/2 & 1/3 & 1/4 & 1/5 & 1/6 \\
1/3 & 1/4 &  1/5 & 1/6 & 1/7 \\
1/4 & 1/5 & 1/6 & 1/7 & 1/8\\
1/5 & 1/6 & 1/7 & 1/8 & 1/9\\
\end{bmatrix}
\end{eqnarray*}

Now, we use \textbf{MatrixPlot} to obtain the image shown in Fig. 1.1.

\begin{figure}

  \includegraphics[width=50mm]{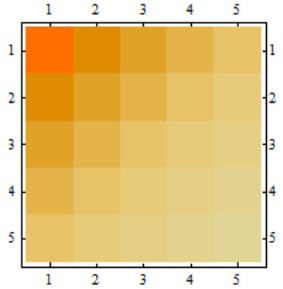}\\
  \caption{2D Plot}\label{Fig. 1}

\end{figure}

\begin{figure}

  \includegraphics[width=70mm]{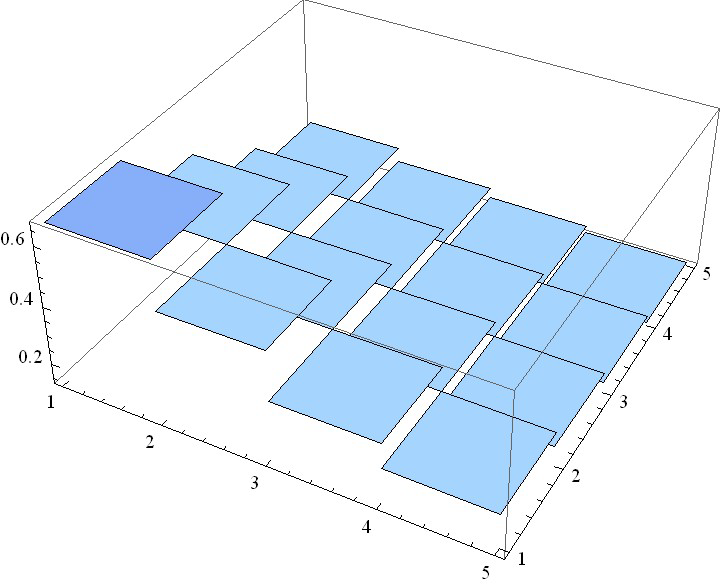}\\
  \caption{3D Plot}\label{Fig. 2}
\end{figure}

With following algorithm which used the Mathematica script\cite{Pet}, we can obtain Hilbert view of a digital image:\\

i. Extract an $n \times n$ subimage $h$ from your favorite greyscale digital image.\\

ii. Multiply subimage $h$ by the corresponding Hilbert $n \times n$ matrix. Let $hilbertim=h.*HilbertMatrix[n]$.\\

iii. Produce a $2D$ MatrixPlot for $hilbertim$ like the one in Fig. 1.1.\\

iv. Use the following scrip to produce a $3D$ plot for $hilbertim$ like the one in Fig. 1.2.

\begin{eqnarray*}
ListPlot3D[HilbertMatrix[n], \rightarrow 0, Mesh \rightarrow None].
\end{eqnarray*}

Again, it can be given cryptography as an example of Hilbert matrix applications. Cryptography is the science of using mathematics to encrypt and decrypt data. Classical cryptanalysis involves an interesting combination of analytical reasoning, application of mathematical tools and pattern finding. In some studies related to cryptographic methods, Hilbert matrix is used for authentication and confidentiality\cite{RajaChak}. It is well known that the Hilbert matrix is very unstable \cite{Matsuo} and so it can be used in security systems.

\section{Hilbert matrix, difference operator and new spaces}

Let $\omega$, $\ell_{\infty}$, $c$, $c_{0}$, $\phi$ denote the sets of all complex, bounded, convergent, null convergent and finite sequences, respectively.
Also, for the sets convergent, bounded and absolutely convergent series, we denote $cs$, $bs$ and $\ell_{1}$.\\

Let $A=(a_{nk})$, $(n,k \in \mathbb{N})$ be an infinite matrix of complex numbers and $X$, $Y$ be subsets of $\omega$. We write $A_{n}x=\sum_{k=0}^{\infty}a_{nk}x_{k}$ and $Ax=(A_{n}x)$ for all $n\in \mathbb{N}$. All the series $A_{n}x$ converge. The set $X_{A}=\{x\in \omega: Ax\in X\}$ is called the matrix domain of $A$ in $X$. We write $(X:Y)$ for the space of those matrices which send the whole of sequence space $X$ into the sequence space $Y$ in this sense.\\

A matrix $A=(a_{nk})$ is called a triangle if $a_{nk}=0$ for $k>n$ and $a_{nn}\neq 0$ for all $n\in \mathbb{N}$. For the triangle matrices $A$, $B$ and a sequence $x$, $A(Bx)=(AB)x$ holds. We remark that a triangle matrix $A$ uniquely has an inverse $A^{-1}=B$ and the matrix $B$ is also a triangle.\\

Let $X$ be a normed sequence space. If $X$ contains a sequence $(b_{n})$ with the property that for every $x\in X$, then there is a unique sequence of scalars $(\alpha_{n})$ such that
\begin{eqnarray*}
\lim_{n\rightarrow\infty}\left\|x-\sum_{k=0}^{n}\alpha_{k}b_{k}\right\|=0.
\end{eqnarray*}
Thus, $(b_{n})$ is called Schauder basis for $X$.\\

It was first used the difference operator in the sequence spaces by K{\i}zmaz\cite{Kizmaz}. The idea of difference sequence spaces of K{\i}zmaz  was generalized by \c{C}olak and Et\cite{ColEt, EtCol}. The difference matrix $\Delta=\delta_{nk}$ defined by

\begin{eqnarray*}
\delta_{nk}:= \left\{ \begin{array}{ccl}
(-1)^{n-k}&, & \quad (n-1\leq k \leq n)\\
0&, & \quad (0<n-1 \textrm{ or } n>k)
\end{array} \right.
\end{eqnarray*}

The difference operator order of $m$ is defined by $\Delta^{m}:\omega \rightarrow \omega$, $(\Delta^{1}x)_{k}=(x_{k}-x_{k-1})$ and $\Delta^{m}x=(\Delta^{1}x)_{k} \circ (\Delta^{m-1}x)_{k}$ for $m\geq 2$.\\

The triangle matrix $\Delta^{(m)}=\delta_{nk}^{(m)}$ defined by
\begin{eqnarray*}
\delta_{nk}:= \left\{ \begin{array}{ccl}
(-1)^{n-k}\binom{m}{n-k}&, & \quad (\max\{0,n-m\}\leq k \leq n)\\
0&, & \quad (0\leq k < \max\{0,n-m\} \textrm{ or } n>k)
\end{array} \right.
\end{eqnarray*}
for all $k,n\in \mathbb{N}$ and for any fixed $m\in \mathbb{N}$.\\

It can be also given Fibonacci matrices as an example of difference matrices\cite{karailk}.\\

The Hilbert matrix is defined by $H_{n}=[h_{ij}]=[\frac{1}{i+j-1}]_{i,j=1}^{n}$ for each $n\in N$. The inverse of Hilbert matrix is defined by
\begin{eqnarray}\label{invHil}
H_{n}^{-1}=(-1)^{i+j}(i+j-1)\binom {n+i-1}{n-j}\binom {n+j-1}{n-i}\binom {i+j-1}{i-1}^{2}
\end{eqnarray}
for all $k,i,j,n\in \mathbb{N}$\cite{Choi}. Polat \cite{Polat}, has defined new spaces by using the Hilbert matrix. Let $h_{c}$, $h_{0}$, $h_{\infty}$ be convergent Hilbert, null convergent Hilbert and bounded Hilbert spaces, respectively. Then, we have:
\begin{eqnarray*}
X=\left\{x\in \omega: Hx\in  Y\right\}
\end{eqnarray*}
where $X=\{h_{c},h_{0},h_{\infty}\}$ and $Y=\{c,c_{0},\ell_{\infty}\}$.\\

Now, we will give new difference Hilbert sequence spaces as below:

\begin{eqnarray*}
&&h_{c}\left(\Delta^{(m)}\right)=\left\{x\in \omega: \Delta^{(m)}x\in  h_{c} \right\}\\
&&h_{0}\left(\Delta^{(m)}\right)=\left\{x\in \omega: \Delta^{(m)}x\in  h_{0} \right\}\\
&&h_{\infty}\left(\Delta^{(m)}\right)=\left\{x\in \omega: \Delta^{(m)}x\in  h_{\infty} \right\}.
\end{eqnarray*}

These new spaces as the set of all sequences whose $\Delta^{(m)}-$transforms of them are in the Hilbert sequences spaces which are defined by Polat\cite{Polat}.\\

We also define the $H\Delta^{(m)}-$transform of a sequence, as below:
\begin{eqnarray}\label{transseq}
y_{n}=\left( H\Delta^{(m)}x \right)_{n}=\sum_{k=1}^{n}\left[\sum_{i=k}^{n}\frac{1}{n+i-1}(-1)^{i-k}\binom{m}{i-k}\right]x_{k}
\end{eqnarray}
for each $m,n\in \mathbb{N}$. Here and after by $H^{(m)}$, we denote the matrix $H^{(m)}=H\Delta^{(m)}=[h_{nk}]$ defined by
\begin{eqnarray*}
h_{nk}=\sum_{i=k}^{n}\frac{1}{n+i-1}(-1)^{i-k}\binom{m}{i-k}
\end{eqnarray*}
for each $k,m,n\in \mathbb{N}$.

\begin{thm}\label{isomorph}
The Hilbert sequences spaces derived by difference operator of $m$ are isomorphic copies of the convergent, null convergent and bounded sequence spaces.
\end{thm}

\begin{proof}
We will only prove that the null convergent Hilbert sequence space is isomorphic copy of the null convergent sequence space. To prove the fact $h_{0}(\Delta^{(m)})\cong c_{0}$, we should show the existence of a linear bijection between the spaces $h_{0}(\Delta^{(m)})$ and $c_{0}$. Consider the transformation $T$ defined, with the notation (\ref{transseq}), from $h_{0}(\Delta^{(m)})$ to $c_{0}$ by $x \rightarrow y=Tx$. The linearity of $T$ is clear. Further, it is trivial that $x=0$ whenever $Tx=0$ and hence $T$ is injective. Let $y\in c_{0}$ and define the sequence $x=(x_{n})$ by
\begin{eqnarray}\label{invseq}
x_{n}=\sum_{k=1}^{n}\left[\sum_{i=k}^{n}\binom{m+n-i-1}{n-i}h_{ik}^{-1}\right]y_{k}
\end{eqnarray}
where $h_{ik}^{-1}$ is defined by (\ref{invHil}). Then,

\begin{eqnarray*}
\lim_{n\rightarrow\infty}\left(H\Delta^{(m)}x\right)_{k}&=&\lim_{n\rightarrow\infty}\sum_{k=1}^{n}\frac{1}{n+k-1}\Delta^{(m)}x_{k}\\
&=&\sum_{k=1}^{n}\frac{1}{n+k-1}\sum_{i=0}^{m}(-1)^{i}\binom{m}{i}x_{k-i}\\
&=&\sum_{k=1}^{n}\left[\sum_{i=k}^{n}\frac{1}{n+i-1}(-1)^{i-k}\binom{m}{i-k}\right]x_{k}=\lim_{n\rightarrow\infty} y_{n}=0.
\end{eqnarray*}
Thus, we have that $x\in h_{0}(\Delta^{(m)}) $. Consequently, $T$ is surjective and is norm preserving. Hence, $T$ is linear bijection which implies that null convergent Hilbert sequence space is isomorphic copy of the null convergent sequence space.\\
\end{proof}

\begin{rem}It is well known that the spaces $c$, $c_{0}$ and $\ell_{\infty}$ are $BK-$spaces. Considering the fact that $\Delta^{(m)}$ is a triangle, we can say that the Hilbert sequences spaces derived by difference operator of $m$, are $BK-$spaces with the norm
\begin{eqnarray}\label{norm}
\|x\|_{\Delta}=\|H\Delta^{(m)}x\|_{\infty}=\sup_{n}\left|\sum_{k=1}^{n}\frac{1}{n+k-1}\sum_{i=0}^{m}(-1)^{i}\binom{m}{i}x_{k-i}\right|.
\end{eqnarray}
\end{rem}

In the theory of matrix domain, it is well known that the matrix domain $X_{A}$ of a normed sequence space $X$ has
a basis if and only if $X$ has a basis whenever $A=(a_{nk})$ is a triangle. Then, we have:

\begin{cor}
Define the sequence $b^{(k)}=(b_{n}^{(k)}(\Delta^{(m)}))_{n\in \mathbb{N}}$ by
\begin{eqnarray*}
b_{n}^{(k)}(\Delta^{(m)}):= \left\{ \begin{array}{ccl}
\sum_{i=k}^{n}\binom{m+n-i-1}{n-i}h_{ik}^{-1}&, & \quad (n\geq k)\\
0&, & \quad (n<k)
\end{array} \right.
\end{eqnarray*}
for every fixed $k\in \mathbb{N}$. The following statements hold:\\

i. The sequence $b^{(k)}(\Delta^{(m)})=(b_{n}^{(k)}(\Delta^{(m)}))_{n\in \mathbb{N}}$ is basis for the null convergent Hilbert sequence space, and for any $x\in h_{0}(\Delta^{(m)})$ has a unique representation of the form
\begin{eqnarray*}
x=\sum_{k}\left(H\Delta^{(m)}x\right)_{k}b^{(k)}.
\end{eqnarray*}\\
ii. The set $\{t, b^{(1)}, b^{(2)},\cdots\}$ is a basis for  the convergent Hilbert sequence space, and for any $x\in h_{c}(\Delta^{(m)})$ has a unique representation of form
\begin{eqnarray*}
x=st+\sum_{k}\left[\left(H\Delta^{(m)}x\right)_{k}-s\right]b^{(k)}
\end{eqnarray*}
where $t=t_{n}(\Delta^{(m)})=\sum_{k=1}^{n}\sum_{i=k}^{n}\binom{m+n-i-1}{n-i}h_{ik}^{-1}$ for all $k\in \mathbb{N}$ and $s=\lim_{k\rightarrow\infty}(H\Delta^{(m)}x)_{k}$.
\end{cor}

If we take into consideration the fact that a space which has a Schauder basis, is separable, then, we can give following corollary:

\begin{cor}
The convergent Hilbert and  null  convergent Hilbert sequence spaces are separable.
\end{cor}

\section{Dual Spaces and Matrix Transformations}

Let $x$ and $y$ be sequences, $X$ and $Y$ be subsets of $\omega$ and $A=(a_{nk})_{n,k=0}^{\infty}$ be an infinite matrix of complex numbers. We write $xy=(x_{k}y_{k})_{k=0}^{\infty}$, $x^{-1}*Y=\{a\in\omega: ax\in Y\}$ and $M(X,Y)=\bigcap_{x\in X}x^{-1}*Y=\{a\in\omega: ax\in Y ~\textrm{ for all }~ x\in X\}$ for the \emph{multiplier space} of $X$ and $Y$. In the special cases of $Y=\{\ell_{1}, cs, bs\}$, we write $x^{\alpha}=x^{-1}*\ell_{1}$,  $x^{\beta}=x^{-1}*cs$,  $x^{\gamma}=x^{-1}*bs$ and $X^{\alpha}=M(X,\ell_{1})$,  $X^{\beta}=M(X,cs)$,  $X^{\gamma}=M(X,bs)$ for the $\alpha-$dual, $\beta-$dual, $\gamma-$dual of $X$. By $A_{n}=(a_{nk})_{k=0}^{\infty}$ we denote the sequence in the $n-$th row of $A$, and we write $A_{n}(x)=\sum_{k=0}^{\infty}a_{nk}x_{k}$ $n=(0,1,...)$ and $A(x)=(A_{n}(x))_{n=0}^{\infty}$, provided $A_{n}\in x^{\beta}$ for all $n$.\\

\begin{lem}\cite[Lemma 5.3]{AB2}\label{duallem1}
Let $X, Y$ be any two sequence spaces. $A\in (X: Y_{T})$ if and only if $TA\in (X:Y)$, where
$A$ an infinite matrix and $T$ a triangle matrix.
\end{lem}

\begin{lem}\label{duallem2}\cite[Theorem 3.1]{AB}
Let $U=(u_{nk})$, be an infinite matrix of complex numbers for all $n,k \in \mathbb{N}$. Let $B^{U}=(b_{nk})$ be defined via a sequence $a=(a_{k})\in\omega$ and inverse of the triangle matrix $U=(u_{nk})$ by
\begin{eqnarray*}
b_{nk}=\sum_{j=k}^na_{j}v_{jk}
\end{eqnarray*}
for all $k,n\in\mathbb{N}$. Then,
\begin{eqnarray*}
X_{U}^{\beta}=\{a=(a_{k})\in\omega: B^{U}\in(X:c)\}.
\end{eqnarray*}
and
\begin{eqnarray*}
X_{U}^{\gamma}=\{a=(a_{k})\in\omega: B^{U}\in(X:\ell_{\infty})\}.
\end{eqnarray*}
\end{lem}

Now, we list the following useful conditions.\\
\begin{eqnarray}\label{eq21}
&&\sup_{n}\sum_{k}|a_{nk}|<\infty \\ \label{eq22}
&&\lim_{n \to \infty}a_{nk}-\alpha_{k}=0\\ \label{eq23}
&&\lim_{n \to \infty}\sum_{k}a_{nk} \quad \textit{exists}\\ \label{eq24}
&&\lim_{n \to \infty}\sum_{k}\left|a_{nk}\right|=\sum_{k}\left|\lim_{n \to \infty}a_{nk}\right| \\ \label{eq27}
&&\lim_{n}a_{nk}=0  \quad \textrm{ for all k }\\ \label{eq11}
&&\sup_{m}\sum_{k}\left|\sum_{n=0}^{m}\right|<\infty \\ \label{eq12}
&&\sum_{n}a_{nk} \quad \textit{convergent for all k}\\ \label{eq13}
&&\sum_{n}\sum_{k}a_{nk} \quad \textit{convergent} \\ \label{eq35}
&&\lim_{n}a_{nk} \quad \textrm{ exists for all k } \\ \label{eq14}
&&\lim_{m}\sum_{k}\left|\sum_{n=m}^{\infty}a_{nk}\right|=0
\end{eqnarray}

\begin{lem}\label{lemMTRX}
For the characterization of the class $(X:Y)$ with
$X=\{c_{0}, c, \ell_{\infty}\}$ and $Y=\{\ell_{\infty}, c, cs, bs\}$, we can give the necessary and sufficient
conditions from Table 1, where

\begin{center}
\begin{tabular}{|l | l | l | l |}
\hline \textbf{1.} (\ref{eq21}) & \textbf{2.} (\ref{eq21}), (\ref{eq35}) & \textbf{3.}  (\ref{eq11}) & \textbf{4.} (\ref{eq11}), (\ref{eq12})  \\
\hline \textbf{5.} (\ref{eq21}), (\ref{eq35}), (\ref{eq23}) &  \textbf{6.} (\ref{eq11}), (\ref{eq12}), (\ref{eq13}) & \textbf{7.} (\ref{eq35}), (\ref{eq24}) & \textbf{8.} (\ref{eq14})\\
\hline
\end{tabular}
\end{center}
\end{lem}

\begin{center}
\begin{tabular}{|c | c c c c|}
\hline
To $\rightarrow$ & $\ell_{\infty}$ & $c$ & bs & cs \\ \hline
From $\downarrow$ &  &  & &\\ \hline
$c_{0}$ & \textbf{1.} & \textbf{2.} & \textbf{3.} & \textbf{4.}\\
$c$ & \textbf{1.} & \textbf{5.} & \textbf{3.} & \textbf{6.}\\
$\ell_{\infty}$ & \textbf{1.} & \textbf{7.} & \textbf{3.} & \textbf{8.}\\
\hline
\end{tabular}

\vspace{0.1cm}Table 1\\

\end{center}

Let $h_{nk}^{-1}$ is defined by (\ref{invHil}). For using in the proof of Theorem \ref{dualthm},
we define the matrix $V=(v_{nk})$ as below:
\begin{eqnarray}\label{mtrxdual}
v_{nk}=\left[\sum_{i=k}^{n}\binom{m+n-i-1}{n-i}h_{ik}^{-1}a_{n}\right]
\end{eqnarray}

\begin{thm}\label{dualthm}
The $\beta-$ and $\gamma-$ duals of the Hilbert sequence spaces derived by difference operator of $m$ defined by
\begin{eqnarray*}
&&\left[h_{c_{0}}\left(\Delta^{(m)}\right)\right]^{\beta}= \left\{ a=(a_{k})\in \omega: V\in (c_{0}:c) \right\}\\
&&\left[h_{c}\left(\Delta^{(m)}\right)\right]^{\beta}= \left\{a=(a_{k})\in \omega: V\in (c:c) \right\}\\
&&\left[h_\infty\left(\Delta{(m)}\right)\right]^{\beta}= \left\{ a=(a_{k})\in \omega: V\in (\ell_{\infty}:c) \right\}\\
&&\left[h_{c_{0}}\left(\Delta^{(m)}\right)\right]^{\gamma}= \left\{a=(a_{k})\in \omega: V\in (c_{0}:\ell_{\infty}) \right\}\\
&&\left[h_{c}\left(\Delta^{(m)}\right)\right]^{\gamma}= \left\{ a=(a_{k})\in \omega: V\in (c:\ell_{\infty}) \right\}\\
&&\left[h_\infty\left(\Delta{(m)}\right)\right]^{\gamma}= \left\{a=(a_{k})\in \omega: V\in (\ell_{\infty}:\ell_{\infty}) \right\}
\end{eqnarray*}
\end{thm}

\begin{proof}
We will only show the $\beta-$ and $\gamma-$ duals of the null convergent Hilbert sequence spaces  derived by difference operator of $m$.
Let $a=(a_{k})\in \omega$. We begin the equality
\begin{eqnarray}\label{dual1}
\sum_{k=1}^{n}a_{k}x_{k}=\sum_{k=1}^{n}\sum_{i=1}^{k}\left[\sum_{j=i}^{k}\binom{m+k-j-1}{k-j}h_{ij}^{-1}\right]a_{k}y_{i}
\end{eqnarray}
\begin{eqnarray*}
&=&\sum_{k=1}^{n}\left\{\sum_{i=1}^{k}\left[\sum_{j=i}^{k}\binom{m+k-j-1}{k-j}h_{ij}^{-1}\right]a_{i}\right\}y_{k}\\
&=&(Vy)_{n}
\end{eqnarray*}
where $V=(v_{nk})$ is defined by (\ref{mtrxdual}). Using (\ref{dual1}), we can see that
$ax=(a_{k}x_{k})\in cs$ or $bs$ whenever $x=(x_{k})\in h_{c_{0}}(\Delta^{(m)})$ if and only if $Vy\in c$ or $\ell_{\infty}$
whenever $y=(y_{k})\in c_{0}$. Then, from Lemma \ref{duallem1} and Lemma \ref{duallem2}, we obtain  the result that
$a=(a_{k})\in \left(h_{c_{0}}(\Delta^{(m)})\right)^{\beta}$ or $a=(a_{k})\in \left(h_{c_{0}}(\Delta^{(m)})\right)^{\gamma}$ if and only if
$V\in (c_{0}:c)$ or $V\in (c_{0}:\ell_{\infty})$, which is what we wished to prove.
\end{proof}

Therefore, the $\beta-$ and $\gamma-$ duals of new spaces will help us in the characterization of the matrix transformations.

Let $X$ and $Y$ be arbitrary subsets of $\omega$. We shall show that,
the characterizations of the classes $(X, Y_{T})$ and $(X_{T},Y)$ can be reduced to that of
$(X, Y)$, where $T$ is a triangle.\\

It is well known that if $h_{c_{0}}(\Delta^{(m)}) \cong c_{0}$, then the equivalence
\begin{eqnarray*}
x\in h_{c_{0}}(\Delta^{(m)}) \Leftrightarrow y\in c_{0}
\end{eqnarray*}
holds. Then, the following theorems will be proved and given
some corollaries which can be obtained to that of Theorems
\ref{mtrxtr1} and \ref{mtrxtr2}. Then, using Lemmas \ref{duallem1} and \ref{duallem2}, we have:

\begin{thm}\label{mtrxtr1}
Consider the infinite matrices $A=(a_{nk})$ and $D=(d_{nk})$. These matrices
get associated with each other the following relations:
\begin{eqnarray}\label{eq1}
d_{nk}= \sum_{j=k}^{\infty}\binom{m+n-j-1}{n-j}h_{jk}^{-1}a_{nj}
\end{eqnarray}
for all $k,m, n\in \mathbb{N}$. Then, the following statements are true:\\
\textbf{i.} $A \in (h_{c_{0}}(\Delta^{(m)}):Y)$ if and only if $\{a_{nk}\}_{k\in\mathbb{N}} \in [h_{c_{0}}(\Delta^{(m)})]^{\beta}$
for all $n\in \mathbb{N}$ and $D\in (c_{0}:Y)$, where $Y$ be any sequence space.\\
\textbf{ii.}  $A \in (h_{c}(\Delta^{(m)}):Y)$ if and only if $\{a_{nk}\}_{k\in\mathbb{N}} \in [h_{c}(\Delta^{(m)})]^{\beta}$
for all $n\in \mathbb{N}$ and $D\in (c:Y)$, where $Y$ be any sequence space.\\
\textbf{iii.}  $A \in (h_{\infty}(\Delta^{(m)}):Y)$ if and only if $\{a_{nk}\}_{k\in\mathbb{N}} \in [h_{\infty}(\Delta^{(m)})]^{\beta}$
for all $n\in \mathbb{N}$ and $D\in (\ell_{\infty}:Y)$, where $Y$ be any sequence space.\\
\end{thm}

\begin{proof}
We assume that the (\ref{eq1}) holds between the entries of $A=(a_{nk})$ and $D=(d_{nk})$.
Let us remember that from Theorem \ref{isomorph}, the spaces $h_{c_{0}}(\Delta^{(m)})$ and $c_{0}$ are linearly isomorphic. Firstly,
we choose any $y=(y_{k})\in c_{0}$ and consider $A \in (h_{c_{0}}(\Delta^{(m)}):Y)$. Then, we obtain that $DH\Delta^{(m)}$ exists and
$\{a_{nk}\}\in \left(h_{c_{0}}\Delta^{(m)}\right)^{\beta}$ for all $k\in \mathbb{N}$. Therefore, the necessity of (\ref{eq1}) yields and
$\{d_{nk}\}\in c_{0}^{\beta}$ for all $k,n\in \mathbb{N}$. Hence, $Dy$ exists for each $y\in c_{0}$. Thus, if we take
$m\rightarrow \infty$ in the equality
\begin{eqnarray*}
\sum_{k=1}^{m}a_{nk}x_{k}=\sum_{k=1}^{m}\left[\sum_{i=1}^{k}\sum_{j=i}^{k}\binom{m+k-j-1}{k-j}h_{ij}^{-1}\right]a_{nk}=\sum_{k}d_{nk}y_{k}
\end{eqnarray*}
for all $m,n\in \mathbb{N}$, then, we understand that $Dy=Ax$. So, we obtain that $D\in (c_{0}:Y)$.\\

Now, we consider that $\{a_{nk}\}_{k\in\mathbb{N}} \in (h_{c_{0}}\Delta^{(m)})^{\beta}$
for all $n\in \mathbb{N}$ and $D\in (c_{0}:Y)$. We take any $x=(x_{k})\in h_{c_{0}}\Delta^{(m)}$. Then, we can see that
$Ax$ exists. Therefore, for $m\rightarrow \infty$, from the equality

\begin{eqnarray*}
\sum_{k=1}^{m}d_{nk}y_{k}=\sum_{k=1}^{m}a_{nk}x_{k}
\end{eqnarray*}
for all $n\in \mathbb{N}$, we obtain that $Ax=Dy$. Therefore, this shows that $A \in (h_{c_{0}}(\Delta^{(m)}:Y)$.
\end{proof}

\begin{thm}\label{mtrxtr2}
Consider that the infinite matrices $A=(a_{nk})$
and $E=(e_{nk})$ with
\begin{eqnarray}\label{eq2}
e_{nk}:=\sum_{k=1}^{n}\sum_{j=k}^{n}\frac{1}{n+j-1}(-1)^{j-k}\binom{m}{j-k}a_{jk}.
\end{eqnarray}
Then, the following statements are true:\\
\textbf{i.} $A=(a_{nk})\in (X:h_{c_{0}}(\Delta^{(m)})$ if and only if $E\in (X :c_{0})$\\
\textbf{ii.} $A=(a_{nk})\in (X:h_{c}(\Delta^{(m)})$ if and only if $E\in (X :c)$\\
\textbf{iii.} $A=(a_{nk})\in (X:h_{\infty}(\Delta^{(m)})$ if and only if $E\in (X :\ell_{\infty})$
\end{thm}

\begin{proof}
We take any $z=(z_{k})\in X$. Using the (\ref{eq2}), we have
\begin{eqnarray}\label{trans1}
\sum_{k=1}^{m}e_{nk}z_{k}=\sum_{k=1}^{m}\left[\sum_{k=1}^{n}\sum_{j=k}^{n}\frac{1}{n+j-1}(-1)^{j-k}\binom{m}{j-k}a_{jk}\right]z_{k}
\end{eqnarray}
for all $m,n\in \mathbb{N}$. Then, for $m\rightarrow \infty$, equation (\ref{trans1}) gives us that $(Ez)_{n}=\{H\Delta^{(m)}(Az)\}_{n}$.
Therefore, one can immediately observe from this that $Az\in h_{c_{0}}(\Delta^{(m)}$ whenever $z\in X$ if and only if $Ez\in c_{0}$ whenever $z\in X$. Thus, the proof is completed.
\end{proof}

\section{Examples}

If we choose any sequence spaces $X$ and $Y$ in Theorem \ref{mtrxtr1} and \ref{mtrxtr2} in previous section,
then, we can find several consequences in every choice. For example, if we take
the space $\ell_{\infty}$ and the spaces which are isomorphic to $\ell_{\infty}$
instead of $Y$ in Theorem \ref{mtrxtr1}, we obtain the following examples:

\begin{exmp}\label{exmpE}
The Euler sequence space $e_{\infty}^{r}$ is defined by $e_{\infty}^{r}=\{x\in \omega: \sup_{n\in\mathbb{N}}|\sum_{k=0}^{n}\binom{n}{k}(1-r)^{n-k}r^{k}x_{k}|<\infty\}$ (\cite{BF2} and \cite{BFM}).
We consider the infinite matrix $A=(a_{nk})$ and define the matrix $C=(c_{nk})$ by
\begin{eqnarray*}
c_{nk}=\sum_{j=0}^{n}\binom{n}{j}(1-r)^{n-j}r^{j}a_{jk}  \quad \quad (k,n\in \mathbb{N}).
\end{eqnarray*}
If we want to get necessary and sufficient conditions for the class $(h_{c_{0}}(\Delta^{(m)}): e_{\infty}^{r})$ in Theorem \ref{mtrxtr1},
then, we replace the entries of the matrix $A$ by those of the matrix $C$.
\end{exmp}

\begin{exmp}\label{exmpR}
Let $T_{n}=\sum_{k=0}^{n}t_{k}$ and $A=(a_{nk})$ be an infinite matrix. We define the matrix $G=(g_{nk})$ by
\begin{eqnarray*}
g_{nk}=\frac{1}{T_{n}}\sum_{j=0}^{n}t_{j}a_{jk}  \quad \quad (k,n\in \mathbb{N}).
\end{eqnarray*}
Then, the necessary and sufficient conditions in order for $A$ belongs to the class $(h_{c_{0}}(\Delta^{(m)}):r_{\infty}^{t})$
are obtained from in Theorem \ref{mtrxtr1} by replacing the entries of the matrix $A$ by those of the matrix $G$;
 where $r_{\infty}^{t}$ is the space of all sequences whose $R^{t}-$transforms is in the space $\ell_{\infty}$ \cite{malk}.
\end{exmp}

\begin{exmp}
In the space $r_{\infty}^{t}$, if we take $t=e$, then, this space become to the Cesaro sequence space of non-absolute type $X_{\infty}$ \cite{NgLee}.
As a special case, Example \ref{exmpR} includes the characterization of class $(h_{c_{0}}(\Delta^{(m)}):r_{\infty}^{t})$.
\end{exmp}

\begin{exmp}
The Taylor sequence space $t_{\infty}^{r}$ is defined by $t_{\infty}^{r}=\{x\in \omega: \sup_{n\in\mathbb{N}}|\sum_{k=n}^{\infty}\binom{k}{n}(1-r)^{n+1}r^{k-n}x_{k}|<\infty\}$ (\cite{kirisci4}).
We consider the infinite matrix $A=(a_{nk})$ and define the matrix $P=(p_{nk})$ by
\begin{eqnarray*}
p_{nk}=\sum_{k=n}^{\infty}\binom{k}{n}(1-r)^{n+1}r^{k-n}a_{jk}  \quad \quad (k,n\in \mathbb{N}).
\end{eqnarray*}
If we want to get necessary and sufficient conditions for the class $(h_{c_{0}}(\Delta^{(m)}): t_{\infty}^{r})$ in Theorem \ref{mtrxtr1},
then, we replace the entries of the matrix $A$ by those of the matrix $P$.
\end{exmp}

If we take the spaces $c$, $cs$ and $bs$ instead of $X$ in Theorem \ref{mtrxtr2}, or $Y$ in Theorem \ref{mtrxtr1}
we can write the following examples. Firstly, we give some conditions and following lemmas:

\begin{eqnarray}\label{eq27x}
&&\lim_{k}a_{nk}=0  \quad \textrm{ for all n }, \\ \label{eq28}
&&\lim_{n \to \infty}\sum_{k}a_{nk}=0, \\ \label{eq26x}
&&\lim_{n \to \infty}\sum_{k} |a_{nk}|=0, \\ \label{eq33}
&&\lim_{n \to \infty}\sum_{k} |a_{nk}-a_{n,k+1}|=0, \\ \label{eq34}
&&\sup_{n}\sum_{k}\left|a_{nk}-a_{n,k+1}\right|<\infty \\ \label{eq36}
&&\lim_{k}\left(a_{nk}-a_{n,k+1}\right) \textrm{ exists for all k } \\ \label{eq37}
&&\lim_{n \to \infty}\sum_{k}\left|a_{nk}-a_{n,k+1}\right|=\sum_{k}\left|\lim_{n \to \infty}(a_{nk}-a_{n,k+1})\right|\\ \label{eq38}
&&\sup_{n}\left|\lim_{k}a_{nk}\right|<\infty
\end{eqnarray}

\begin{lem}
Consider that the $X\in\{\ell_{\infty}, c, bs, cs\}$ and $Y\in \{c_{0}\}$.
The necessary and sufficient conditions for $A\in (X:Y)$ can be read the following, from Table 2:

\begin{center}
\begin{tabular}{|l | l | l | l |}
\hline \textbf{9.} (\ref{eq26x}) & \textbf{10.} (\ref{eq21}), (\ref{eq27}), (\ref{eq28}) & \textbf{11.} (\ref{eq27x}), (\ref{eq33}) & \textbf{12.} (\ref{eq27}), (\ref{eq34})  \\
\hline \textbf{13.} (\ref{eq27x}), (\ref{eq36}), (\ref{eq37}) &  \textbf{14.} (\ref{eq34}), (\ref{eq35})  & \textbf{15.} (\ref{eq27x}), (\ref{eq34}) & \textbf{16.} (\ref{eq34}), (\ref{eq38}) \\ \hline
\end{tabular}
\end{center}

\end{lem}

\begin{center}
\begin{tabular}{|c | c c  c c|}
\hline
From $\rightarrow$ & $\ell_{\infty}$ & $c$ &  $bs$ & $cs$ \\ \hline
To $\downarrow$ &    & & &\\ \hline
$c_{0}$ & \textbf{9.} & \textbf{10.} & \textbf{11.} & \textbf{12.}\\
$c$ & \textbf{7.} & \textbf{5.} & \textbf{13.} & \textbf{14.}\\
$\ell_{\infty}$ & \textbf{1.} & \textbf{1.} & \textbf{15.} & \textbf{16.}\\
\hline
\end{tabular}

\vspace{0.1cm}Table 2\\
\end{center}

\begin{exmp}
We choose $X\in \{h_{c_{0}}(\Delta^{(m)})\}, h_{c}(\Delta^{(m)}), h_{\ell_{\infty}}(\Delta^{(m)})$ and $Y\in \{\ell_{\infty}, c, cs, bs, f\}$.
The necessary and sufficient conditions for $A\in (X:Y)$ can be taken from the Table 3:
\end{exmp}
\begin{itemize}
  \item[\textbf{1a.}] (\ref{eq21}) holds with $d_{nk}$ instead of $a_{nk}$.
 \item[\textbf{2a.}] (\ref{eq21}), (\ref{eq35}) hold with $d_{nk}$ instead of $a_{nk}$.
\item[\textbf{3a.}] (\ref{eq11}) holds with $d_{nk}$ instead of $a_{nk}$.
\item[\textbf{4a.}] (\ref{eq11}), (\ref{eq12}) hold with $d_{nk}$ instead of $a_{nk}$.
\item[\textbf{5a.}] (\ref{eq21}), (\ref{eq35}), (\ref{eq23}) hold with $d_{nk}$ instead of $a_{nk}$.
\item[\textbf{6a.}] (\ref{eq11}), (\ref{eq12}), (\ref{eq13}) hold with $d_{nk}$ instead of $a_{nk}$.
\item[\textbf{7a.}] (\ref{eq35}), (\ref{eq24}) hold with $d_{nk}$ instead of $a_{nk}$.
\item[\textbf{8a.}] (\ref{eq14}) holds with $d_{nk}$ instead of $a_{nk}$.
\end{itemize}

\begin{center}
\begin{tabular}{|c | c c c c |}
\hline
To $\rightarrow$ & $\ell_{\infty}$ & $c$ & $bs$ & $cs$ \\ \hline
From $\downarrow$ &  &  & & \\ \hline
$h_{c_{0}}(\Delta^{(m)})$ & \textbf{1a.} & \textbf{2a.} & \textbf{3a.} & \textbf{4a.} \\
$h_{c}(\Delta^{(m)})$ & \textbf{1a.} & \textbf{5a.} & \textbf{3a.} & \textbf{6a.} \\
$h_{\ell_{\infty}}(\Delta^{(m)})$ & \textbf{1a.} & \textbf{7a.} & \textbf{3a.} & \textbf{8a.}\\
\hline
\end{tabular}

\vspace{0.1cm}Table 3\\
\end{center}

\begin{exmp}
Consider that the $X\in\{\ell_{\infty}, c, bs, cs\}$ and $Y\in \{h_{c_{0}}(\Delta^{(m)})\}, h_{c}(\Delta^{(m)}), h_{\ell_{\infty}}(\Delta^{(m)})$.
The necessary and sufficient conditions for $A\in (X:Y)$ can be read the following, from Table 4:
\begin{itemize}
\item[\textbf{9a.}] (\ref{eq26x}) holds with $e_{nk}$ instead of $a_{nk}$.
\item[\textbf{10a.}] (\ref{eq21}), (\ref{eq27}), (\ref{eq28}) hold with $e_{nk}$ instead of $a_{nk}$.
\item[\textbf{11a.}] (\ref{eq27x}), (\ref{eq33}) hold with $e_{nk}$ instead of $a_{nk}$.
\item[\textbf{12a.}] (\ref{eq27}), (\ref{eq34}) hold with $e_{nk}$ instead of $a_{nk}$.
\item[\textbf{13a.}] (\ref{eq27x}), (\ref{eq36}), (\ref{eq37}) hold with $e_{nk}$ instead of $a_{nk}$.
\item[\textbf{14a.}] (\ref{eq34}), (\ref{eq35}) hold with $e_{nk}$ instead of $a_{nk}$.
\item[\textbf{15a.}] (\ref{eq27x}), (\ref{eq34}) hold with $e_{nk}$ instead of $a_{nk}$.
\item[\textbf{16a.}] (\ref{eq34}), (\ref{eq38}) hold with $e_{nk}$ instead of $a_{nk}$.
\end{itemize}
\end{exmp}

\begin{center}
\begin{tabular}{|c | c c  c c|}
\hline
From $\rightarrow$ & $\ell_{\infty}$ & $c$ &  $bs$ & $cs$ \\ \hline
To $\downarrow$ &    & & &\\ \hline
$h_{c_{0}}(\Delta^{(m)})$ & \textbf{9a.} & \textbf{10a.} & \textbf{11a.} & \textbf{12a.}\\
$h_{c}(\Delta^{(m)})$ & \textbf{7a.} & \textbf{5a.} & \textbf{13a.} & \textbf{14a.} \\
$h_{\ell_{\infty}}(\Delta^{(m)})$ & \textbf{1a.} & \textbf{1a.} & \textbf{15a.} & \textbf{16a.}\\
\hline
\end{tabular}

\vspace{0.1cm}Table 4\\
\end{center}

\section{Conclusion}
In 1894, Hilbert introduced the Hilbert matrix. The Hilbert matrices are notable examples of poorly conditioned
(ill-conditioned) matrices, making them notoriously difficult to use in numerical computation. That is, Hilbert matrices whose entries are
specified as machine-precision numbers are difficult to invert using numerical techniques. That's why, we offer some examples related to the using Hilbert matrix such as image processing and cryptography. We also give an algorithm. Further, we construct new spaces with the Hilbert matrix and difference operator of order $m$. We calculate dual spaces of new spaces and characterize some matrix classes. In last section, we give some examples of matrix classes. Images of new spaces can be plotted using the Mathematica as a continuation of this study. Again, different applications of cryptography can be investigated.

\section*{Conflict of Interests}
The authors declare that there are no conflict of interests regarding the publication of this paper.

\end{document}